\documentclass[final]{siamltex}
\usepackage[hmargin=1in,vmargin=1.2in]{geometry}
\usepackage{graphicx}
\usepackage{float}
\usepackage{caption}
\usepackage{subfigure}
\usepackage{amsfonts,amssymb}
\usepackage{color}
\definecolor{myblue}{rgb}{0.2, 0.5, 0.8}

\usepackage[pagewise]{lineno}
\usepackage{color}
\definecolor{darkgreen}{rgb}{0.1,0.5,0.1}
\definecolor{darkblue}{rgb}{0.2,0.2,1.0}
\usepackage[colorlinks=true,linkcolor=darkblue,citecolor=darkblue,
            filecolor=darkblue,urlcolor=darkgreen]{hyperref}

\usepackage{amsmath}

\newcommand{\p}{\partial}
\newcommand{\fr}{\frac}

\title{Computational study of shock waves propagating through air-plastic-water interfaces.}
\author{Mauricio J. Del Razo\thanks{Department of Applied Mathematics, 
        University of Washington, Seattle, WA 98195-3925 ({\tt maojrs@uw.edu, rjl@uw.edu}).} \and Randall J. LeVeque \footnotemark[1]}
        
\begin{document}

% \date{}

\maketitle

\begin{abstract}
 The following study is motivated by experimental studies in traumatic brain injury (TBI). Recent research has demonstrated that low intensity non-impact blast wave exposure frequently leads to mild traumatic brain injury (mTBI); however, the mechanisms connecting the blast waves and the mTBI remain unclear. Collaborators at the Seattle VA Hospital are doing experiments to understand how blast waves can produce mTBI. In order to gain insight that is hard to obtain by experimental means, we have developed conservative finite volume methods for interface-shock wave interaction to simulate these experiments. A 1D model of their experimental setup has been implemented using Euler equations for compressible fluids. These equations are coupled with a Tammann equation of state (EOS) that allows us to model compressible gas along with almost incompressible fluids or elastic solids. A hybrid HLLC-exact Eulerian-Lagrangian Riemann solver for Tammann EOS with a jump in the parameters has been developed. The model has shown that if the plastic interface is very thin, it can be neglected. This result might be very helpful to model more complicated setups in higher dimensions. 
\end{abstract}



\pagestyle{myheadings} \markboth{Computational study of shock waves propagating through air-plastic-water interfaces}{M. J. Del Razo \& R. J. LeVeque} 

% % % % % % % % % % % % % % % % % % % % % % % % % % % % % % % % % % % % % % % % % % % % % % % % % % % % % % % % % % % % % % % % % % % % % % % % % % % % 
\section{Introduction} 
Recent research has demonstrated that low intensity non-impact blast wave exposure frequently leads to mild traumatic brain injury (mTBI)\cite{Peskind2011,xiong2013animal}; however, the mechanisms connecting the blast waves and the mTBI remain unclear. Conventional imaging techniques, like computerized tomography and MRI, fail to reveal the damage when the injuries are mild. This suggests that the injury mechanisms might occur at very small length scales, even at the scale of a single cell. Some small scale damage could include metabolic cells' changes or even alterations in the blood brain barrier(BBB) diffusion, as hinted in \cite{Kucherov2012}. 

\begin{figure}[H]
\centering
\includegraphics[width=0.8\textwidth]{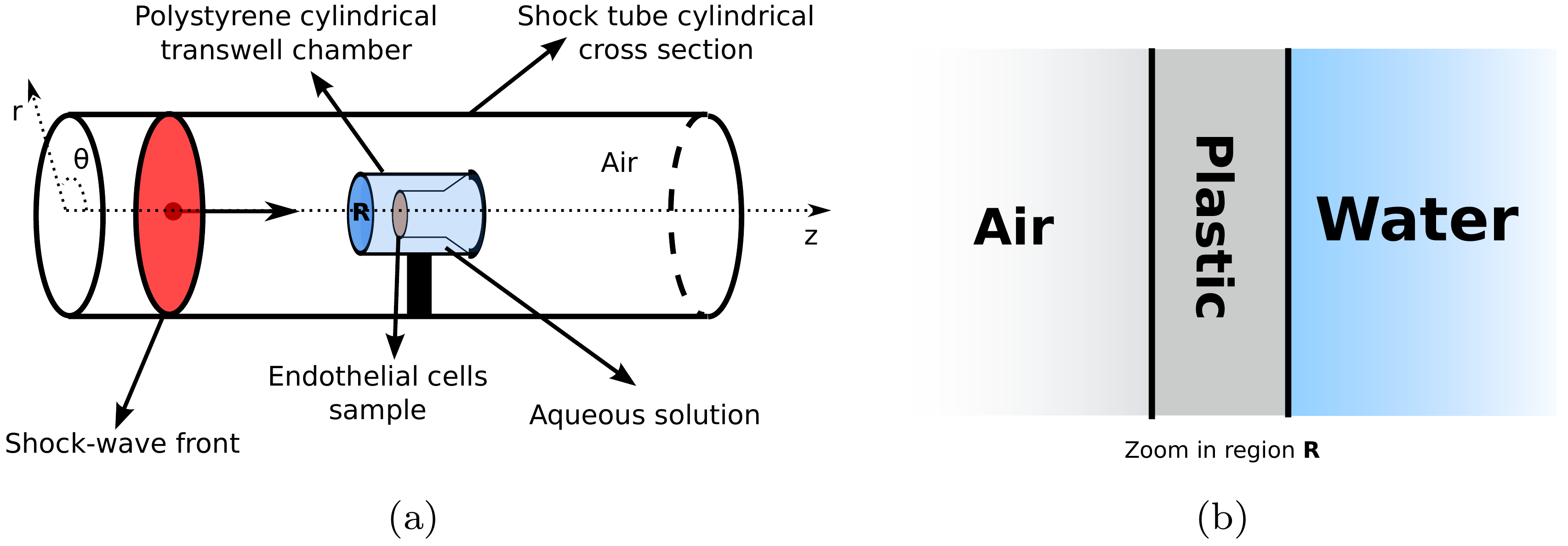} 
\caption{(a) Cartoon of experimental setup. A cylindrical cross section of the shock tube with the sample
is shown. The shock wave fronts travels through air in the shock tube hitting the plastic first, then the water
and finally the endothelial cells sample. (b) Zoom-in of the Region \textbf{R} of (a). This region is the main concern
of this work, and it will be modeled as one dimensional.} 
\label{fig:stube}
\end{figure}

In order to gain insight on the injury mechanisms, experimentalists at the Seattle Veterans Administration hospital (VA hospital) are trying to reproduce the effects of blast waves in brain cells with the aid of a shock-tube. A sample with endothelial cells, the cells that form the BBB, is prepared and placed inside the shock-tube. The differences in several properties of the cells are studied before and after the blast wave has passed. The endothelial cells sample is placed inside of a polystyrene (plastic) cylindrical transwell container filled with an aqueous solution; a cartoon of the experimental setup  is presented in Figure \ref{fig:stube}a. 

The shock-tube is a long cavity filled with air, where the shock wave pressure profile can be measured. However, we don't know what will happen when the shock wave interacts with the polystyrene cylindrical transwell filled with water. The physical effects of the shock wave crossing from air to plastic and from plastic to water are not evident from experimental data nor easy to obtain through experimental techniques. The model presented here
concentrates on the one dimensional idealized scenario (Figure \ref{fig:stube}b) of the region \textbf{R} of Figure \ref{fig:stube}a. The computational simulations presented will help bring insight on the behavior of the shock wave crossing air-plastic-water interfaces, as the one 
shown in Figure \ref{fig:stube}b. They will also provide quantitative data of the strength of the shock wave before and after crossing the interfaces, and they will show that if the plastic interface is thin enough, it can be completely neglected.

In Section \ref{sec:exactac}, we show an exact result for linear acoustics of how the plastic interface becomes irrelevant as its thickness goes to zero. In Section \ref{sec:TBInumimp}, the appropriate model to represent the experimental setup will be presented; an overview of its numerical implementation is presented in Section \ref{sec:nummethods}. In Section \ref{sec:TBI:compexp}, two computational experiments are described to show the plastic interface can be neglected without losing much accuracy. Moreover, we gained insights in the behavior of the shock wave when crossing the interfaces, of which the most relevant is the increase in pressure when the shock wave crosses from air to water. A discussion of the results is given on Section \ref{sec:TBIdisc}.
These results form a basis for on-going work on this application, to be described
further in future publications \cite{delrazo2015_02,delrazo2015_01}.

\section{Air-Plastic-Water interface for linear acoustics}
\label{sec:exactac}
Before trying to model the shock waves with complicated non-linear systems of conservation laws, we can start with a very simple case to gain insight: linear acoustics. In this case, we can compute the exact solution of the transmitted pressure
through the air-plastic-water interface as a function of the acoustic impedance of each material and the plastic width.
This can be derived from the fact that an acoustic wave with incident pressure $p_0$ on the left of an interface between medium $A$ (left) and $B$(right) produces a reflected and a transmitted wave with pressures given by
\begin{align*}
  p_T=p_0\frac{2Z_B}{Z_A + Z_B} \hspace{10mm} p_R=p_0\frac{Z_B-Z_A}{Z_A+Z_B},
\end{align*}
where $Z_k$ denotes the acoustic impedance of medium $k$. These relations can be easily derived from linear acoustics \cite{randysrbook}.
Now consider a triple one dimensional interface: air-plastic-water. With this
setup, there will be an infinite number of reflections in the plastic layer. 
The transmitted wave in water after the $N^{th}$ reflection through the plastic yields
\begin{align*}
  p_T^{N^{th}} = \frac{2Z_w}{Z_w+Z_p}\left(\frac{Z_a-Z_p}{Z_a+Z_p}\right)^{N-1} 
                                     \left(\frac{Z_w-Z_p}{Z_w+Z_p}\right)^{N-1}
                 \frac{2Z_p}{Z_p+Z_a} p_0,
\end{align*}
where $Z_a,Z_p,Z_w$ are the air, plastic and water impedances. 
Each transmitted wave increases the pressure behind the initially transmitted wave
slightly and the asymptotic final amplitude of the transmitted
wave is given by the sum of all these contributions, i.e:
\begin{gather*}
  p_{T_{water}} = \sum_{N=1}^{\infty} p_T^{N^{th}} = 
\frac{4Z_w Z_p p_0}{(Z_w+Z_p)(Z_p+z_a)}\sum_{N=0}^{\infty} 
   \left(\frac{(Z_a-Z_p)(Z_w-Z_p)}{(Z_a+Z_p)(Z_w+Z_p)}\right)^{N}.
\end{gather*}
Summing this geometric series yields
\begin{align*}
  p_{T_{water}} = p_0\frac{2Z_w}{Z_w + Z_a}.
\end{align*}
When the plastic layer is very thin this asymptotic value is quickly reached,
and we note that it is
exactly the same as if the plastic interface didn't exist and we computed the
transmission coefficient directly from air into water. 
Note we assumed the pressure profile on the left was a constant, $p_0$. However, this can be more complicated. It can have a decaying tail, in which case there will be interference from the tail in the reflected and transmitted waves. Nonetheless, assuming the plastic width is $w_0$, the time elapsed between two transmitted waves in the water interface is given by $\tau=2w_0/c_p$, where $c_p$ is the speed sound in plastic. Therefore, as $w_0\rightarrow 0$, the elapsed time does too. As a consequence, the interference from
the tail will also disappear and the plastic interface can be neglected without losing accuracy.

This calculation shows that for linear acoustics, if the plastic interface is very thin in comparison to the experiment characteristic length scales, the plastic interface can be neglected without significant loss of accuracy. From this result, we can expect similar results for the one-dimensional non-linear system that we will employ to model the experimental setup.

\section{The model}
\label{sec:TBInumimp}

In order to perform a computational experiment, we first need to determine which are the most appropriate mathematical
equations to describe the experimental setup. As we are interested in pressure waves (sound-waves), a
first approach is to employ linear acoustic equations. However, as the form of the equations is linear,
shock wave formation is not possible. A viable alternative is the Euler equations for 
compressible inviscid flow. As the equations are nonlinear, shock wave formation is accurately modeled. 
Furthermore, Euler equations also model conservation of energy, whose connection with temperature might be relevant to cell injury. Additionally, for this experimental setup, we are not concerned with large-scale movement of the fluid, so viscosity can be neglected by employing the inviscid equations. One more issue is that air has a high compressibility while water and plastic are almost incompressible. Using different parameters for the equations of state (EOS) for each material, we can model the three materials with the same equations. The equations are solved using the methods briefly explained in Section \ref{sec:nummethods}.

A first approach to this problem is to use one-dimensional Euler equations by focusing in a cross section close to 
the central axis of Figure \ref{fig:stube},
\begin{gather}
\label{eq:Eulercyl}
\begin{gathered}
  \fr{\p}{\p t}
    \left[\begin{array}{c} \rho \\   \rho u \\    E   \end{array} \right] 
+ \fr{\p}{\p x} 
    \left[\begin{array}{c} \rho u \\   \rho u^2 + p \\   u(E+p)   \end{array} \right] 
=  0,
\end{gathered}
\end{gather}
where $\rho$ is the density, $u$ denote the velocities in the $x$ (shock tube axis) direction, $p$ is the pressure and $E$ the internal energy.

\subsection{Tammann equations of state}
\label{sec:EOS}
The system of equations (\ref{eq:Eulercyl}) is closed with the addition of an EOS. It is usually given as a relation 
between pressure, density and specific internal energy, i.e. $p=p(\rho,e)$.  The most well known EOS is the 
one for an ideal gas $p = (\gamma - 1) \rho e$, where $\gamma$ is the heat capacity ratio, $\rho$ is 
the density and $e$ the specific internal energy.  While this 
EOS is very good to describe the behavior of gasses, it might not be appropriate 
to model nearly incompressible materials like water or elastic solids.

Several alternatives exist. In this work, we will use the stiffened gas EOS (SGEOS) or also 
known as Tammann EOS. This equation of state is very useful to model a wide range of fluids even in the 
presence of strong shock waves \cite{Fagnan2012}. The Tammann EOS is given by
\begin{align}
  p = (\gamma - 1) \rho e -\gamma p_{\infty},
  \label{eq:SGEOS}
\end{align}
where $\gamma$ and $p_{\infty}$ can be determined experimentally for different materials. The Tammann EOS
and the ideal gas EOS are the same except for the extra term $-\gamma p_{\infty}$. As the pressure is a positive quantity, 
and $\gamma,p_\infty > 0$, this EOS can only be aimed for higher density fluids. For fluids with $p_\infty \gg 1$,
the relative change in density, when changing the pressure, is very small. Consequently, the Tammann EOS
is a good approximation for fluids with a very low compressibility. It can also be helpful to 
model elastic solids, like plastic. It's worth mentioning that for sufficiently weak shocks the Tammann EOS can be further 
simplified to the Tait EOS, see \cite{Fagnan2012}.

\section{Numerical implementation}
\label{sec:nummethods}

The Euler equations are a hyperbolic system of conservation laws, so they can be solved employing finite volume methods (FVM). This is done by employing the wave propagation algorithms described in \cite{randysrbook} and implemented in Clawpack \cite{clawpack}. The fundamental problem to solve at each cell interface of our computation is the well known Riemann problem. The equations of motion are solved by implementing a hybrid Riemann HLLC-exact type approximate solver for one-dimensional Euler equations with interfaces. This solver couples a Eulerian HLLC approximate Riemann solver to an exact Riemann solver for the Tammannn EOS in Lagrangian coordinates. As the interfaces are represented by contact discontinuities, the HLLC solver is ideal to deal accurately with interface problems. The method can easily be extended to two dimensions by employing dimensional splitting. An overview of the HLLC Riemann solver and the coupling is given next.

\subsection{The HLLC solver}
A general one-dimensional Riemann problem 
for a system of conservation laws like Euler equations can be stated as
\begin{align}
  \mathbf{Q} + \mathbf{F}(\mathbf{Q})_x = 0, \ \ \ \ \ \
  \mathbf{Q}(x,0) = \left\{ 
  \begin{array}{l l}
   \mathbf{Q}_L  \  \text{if $x<0$} \\
   \mathbf{Q}_R  \  \text{if $x>0$}, 
  \end{array} \right. \label{eq:rphll}
\end{align}
where $\mathbf{Q}$ is the vector of state variables.
The HLLC (Harten-Lax-van Leer-Contact) solver is an approximate Riemann solver to the problem (\ref{eq:rphll}).
As the Riemann solution to the one-dimensional Euler equations consists of three
waves: two acoustic waves with a 
contact discontinuity in between.  The main idea of the HLLC solver is, given
upper bounds for the left and right going wave speeds $S_L$ and $S_R$, 
assume a wave configuration of three waves separating four constant states. The approximate solution for this method will be of the form
\begin{align*}
  \tilde{\mathbf{Q}}&(x,t) = \left\{ 
  \begin{array}{l l}
   \mathbf{Q}_L & \  \text{if $\fr{x}{t} \le S_L$},\\
   \mathbf{Q}_{*L} & \  \text{if $S_L \le \fr{x}{t} \le S_*$},\\
   \mathbf{Q}_{*R} & \  \text{if $S_* \le \fr{x}{t} \le S_R$},\\
   \mathbf{Q}_R & \  \text{if $\fr{x}{t}\ge S_R$},
  \end{array} \right.
\end{align*}
where $S_*$ is the approximate wave speed of the contact discontinuity. Assuming we can obtain $S_L$ and $S_R$ by some other algorithm, we 
only need to find $\mathbf{Q}_{*L}$, $\mathbf{Q}_{*R}$ and $S_*$ to solve the problem. These quantities can be obtained
by integrating over a box in the $x,t$ plane, using the Rankine-Hugoniot conditions and assuming constant pressure and normal velocity across the contact discontinuity, see \cite{torosbook}. The desired states 
for the Euler equations with $\mathbf{Q} = [\rho,\rho u, E]^T$ are given by
\begin{align}
  \mathbf{Q}_{*k} = \fr{S_k \mathbf{Q}_k - \mathbf{F}_k + p_* \mathbf{D}}{S_L - S_*}, \ \ \
  \text{with:} \ \ \ D = [0,1,S_*]^T ,
  \label{eq:solRP}
\end{align}
 where the contact discontinuity speed is given by
\begin{align*}
  S_* = \fr{p_R - p_L + \rho_L u_L (S_L - u_L) -\rho_R u_R(S_R - u_R)}{\rho_L(S_L - u_L) - \rho_R(S_R - u_R)},
\end{align*}
and the pressure in the middle state $p_*$ is
\begin{align*}
 p_* = p_L + \rho_L(S_L - u_L)(S_* - u_L) = p_R + \rho_R(S_R-u_R)(S_*-u_R)
\end{align*}

where $\rho_k$, $u_k$ with $k=L,R$ are the left or right density and speed in the Euler equations, for detailed calculations the reader is referred to \cite{torosbook}.

In order to calculate the wave speeds $S_L$ and $S_R$, we will
need to calculate the sound speed. This is where we require the EOS. A simple estimate is the one 
given by Davis \cite{torosbook} as 
$S_L = \text{min}\{u_L - c_L,u_R - c_R\}$ and $S_R = \text{max}\{u_L + c_L, u_R + c_R\}$,
where $u_k$ is the normal velocity and $c_k$ is the sound speed on each side, $k=L,R$. The
speed of sound using the Tammann EOS is given by $c = \sqrt{\gamma\fr{p + p_{\infty}}{\rho}}$.
A possible improvement is to employ Roe averages in wave speed estimates
\cite{einfeldt1988godunov}. These Roe averages
can be calculated using different configurations, some might be more accurate when dealing with 
interfaces, as pointed out in \cite{Hu2009}

\begin{figure}[H]
\centering
\includegraphics[width=0.7\textwidth]{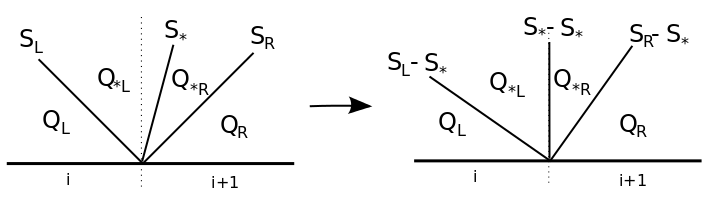}
\caption{Transformation from Eulerian coordinates to Lagrangian coordinates for the HLLC or exact Riemann solver between grid cells $i$ and $i+1$. The transformation is only done when one or more of the cells are part of the highly incompressible
material where $S_* \ll 1$.} 
\label{fig:hllc}
\end{figure}

The HLLC solver just discussed works well for solving the one dimensional Euler equations with an ideal gas EOS. However, we want to implement the HLLC solver with the Tammann EOS across an air-plastic-water interface. The difference between the parameters for different materials in the Tammann EOS is of several orders of magnitude, as shown in Table \ref{tab:param}. This generates instabilities in the HLLC solver, more so in a two or higher dimensional setting. The instability is due to the fact that we model the interfaces as being fixed in space; however, there is always a displacement of the contact discontinuity, i.e. the interface. In order to solve this issue, we model the air in Eulerian coordinates using the usual HLLC solver. If any of the cells is a water or plastic cell, we modify our original HLLC or exact solver to work in Lagrangian coordinates, where the interface is actually fixed with respect to the reference frame. This is done by displacing the frame of reference by $S_*$, 
\begin{align}
 \tilde{S_L} = S_L - S_* \ \ \ \ \ \ \
 \tilde{S_*} = 0 \ \ \ \ \ \ \
 \tilde{S_R} = S_R -  S_* 
 \label{eultolag}
\end{align}
For instance, assume we are running a one dimensional simulation of the Euler
equations, with a fixed interface modeled by a jump in the parameters of the EOS. The interface is aligned to the edge between cells $i$  and $i+1$, the transformed Riemann solver will be as shown in Figure \ref{fig:hllc}. This will ensure the contact discontinuity velocity is zero and consequently the interface is modeled as fixed. The wave contributions will be the correct ones, since we are just modifying the wave velocity and not the solution $q$'s. There is, of course, an error made at the interface when coupling the two descriptions; however, as the displacements of the interface are very small due to very low compressibility, this error is not big enough to cause instabilities as before. The speeds (\ref{eultolag}) along with the solution in the different regions (\ref{eq:solRP}) are the output of our Riemann solver into the Finite Volume method implemented into Clawpack \cite{clawpack}. Note the transformed speeds in (\ref{eultolag}) should only be used in the presence of an interface; otherwise, the Riemann solver should output 
the non-transformed speeds $S_L, S_*$ and $S_R$.

In order to provide better accuracy in certain regions, we also developed and implemented an exact Riemann solver for the
Euler equations using the Tammann EOS with a jump in the parameters, similar to the one by Ivings \& Toro \cite{Ivings1998}.
For the sake of brevity, this will be presented in a forthcoming manuscript. The HLLC solver will be used to model the air in Eulerian coordinates, and the exact solver will be used to model the interface and the almost incompressible material in Lagrangian coordinates. The transformation to Lagrangian coordinates for the exact solver is exactly the same one as in equations (\ref{eultolag}).

\section{Computational experiments}
\label{sec:TBI:compexp}
We want to study the importance of the plastic interface in the computational experiment 
illustrated in Figure \ref{fig:stube}. The simplest scenario is a one-dimensional model of the 
Region \textbf{R} of the transwell in Figure \ref{fig:stube}b, i.e an air-plastic-water interface. 
In this section, we will study a shock wave crossing a general air-plastic-water interface 
for different plastic widths. This will bring insight on the behavior of the shock wave, and it 
will show how relevant is the thin plastic interface. We will analyze the one-dimensional 
air-plastic-water interface problem, and then we'll compare it to the simpler one-dimensional 
air-water interface.  

\subsection{Air-plastic-water interface}
\begin{figure*}[t]
  \centering
  \includegraphics[width=0.9\textwidth]{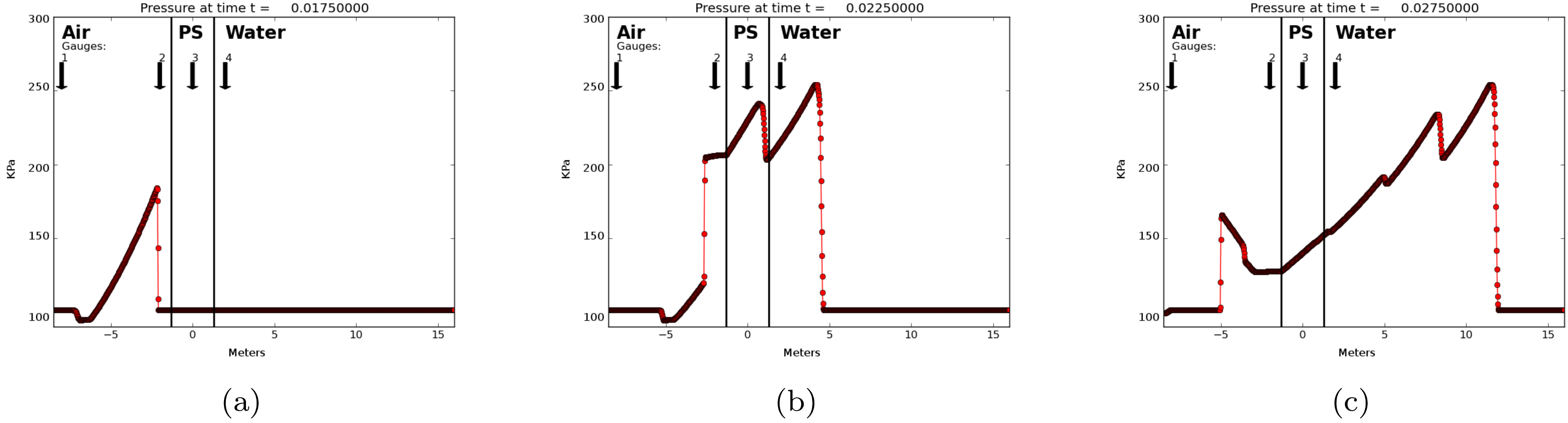} 
  \caption{Shock wave crossing the air-plastic-water interface at different times. The arrows indicate the position
of the 4 gauges that measure the pressure as a function of time. The gauges are numbered from left to right, and the plastic interface width for this case is $2.6$m }
  \label{fig:TBIanim}
\end{figure*}

\begin{figure}[H]
\centering
\includegraphics[width=0.38\paperwidth]{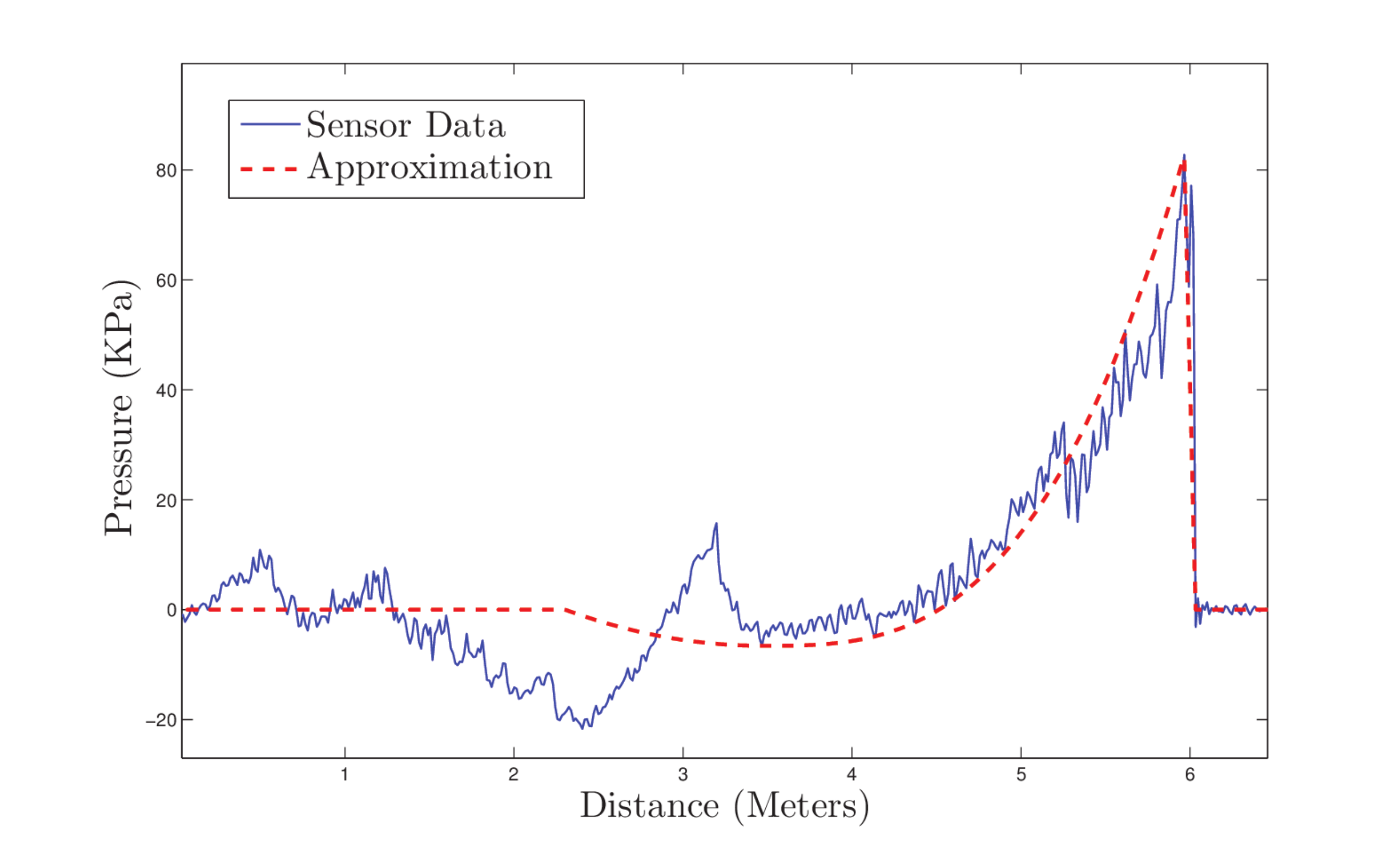}
\caption{The shock wave form obtained from the sensor inside the shock tube
is shown as the solid thin line. The coarse approximation 
to be used as an initial condition in our simulation is shown with a dashed line. An 
average speed of sound of $c=344$ $m$/$s$ is assumed.} 
\label{fig:dataTBI}
\end{figure}

The first step is to input the right initial conditions into our simulation. The actual form of the shock wave traveling through
the shock tube before hitting the transwell was obtained experimentally. The
sensor outputs pressure amplitude as a function of time. Assuming an average
speed of sound in air, it can be converted to a function of distance as shown in
Figure \ref{fig:dataTBI}. The shape can be broadly approximated by an idealized
shock wave (dashed line in Figure \ref{fig:dataTBI}).  This approximated form of
the shock wave is introduced as part of the initial condition on the left part of the domain 
in the simulation; however,
this is not a trivial procedure since we must input the density, momentum and energy, 
and we only have the pressure. Using the isentropic EOS, the ideal gas EOS and 
the expression for the speed of sound, an educated guess for the initial condition 
in terms of the pressure is given far away from the transwell, elsewhere the 
initial condition will be constant: constant density, zero velocity and atmospheric 
pressure. This initial 
condition is then modified until the right amplitude and shape of the shock wave 
front is obtained. The resulting shape of the shock wave before hitting the 
interface can be seen in Figure \ref{fig:TBIanim}a. The pressure is measured in 
KPa with an ambient base pressure of $1\text{ATM}=101.325\text{KPa}$.

The one-dimensional equations with the triple interface are solved using the methods 
mentioned in Section \ref{sec:TBInumimp}. The different materials are modeled using 
different parameters for the Tammann EOS, see Section \ref{sec:EOS}. The choice of 
parameters is shown in Table \ref{tab:param}.
\begin{table}[H]
  \centering
  \begin{tabular}{| l || c | c |}
    \hline
     Material                & $\gamma$     & $p_\infty (GPa)$ \\ \hline 
     Air (Ideal gas EOS)     & 1.4          & 0.0             \\ \hline
     Plastic (polystyrene)   & 1.1          & 4.79            \\ \hline
     Water                   & 7.15         & 0.3             \\
    \hline
  \end{tabular}
  \caption{Parameters for the Tammann EOS to model the different materials. 
  The parameters for air and water were taken from \cite{Fagnan2012}. Since the polystyrene
  is a solid, $\gamma$ was chosen to be close to 1, and $p_\infty$ was adjusted to yield 
  the right speed of sound in polystyrene.}
  \label{tab:param}
\end{table}

The solution of the shock wave crossing the three interfaces; air, plastic and water, 
at different times is shown on Figure \ref{fig:TBIanim}. It can be observed that every time the pressure
wave hits an interface, part of the wave is reflected and part of it is transmitted. This effect can occur multiple times depending on how the interfaces are set up. We can also observe that 
the amplitude of the shock wave increases as it passes from air to plastic and decreases 
when passing from plastic to water. This effect is mostly due to the continuity of pressure at the change in compressibility. In order to keep the pressure at the interface 
continuous, the transmitted wave amplitude has to be the same as the sum of the incident wave and 
the reflected wave. When the compressibility is very high in the adjacent material, the interface 
will behave similarly to a solid wall. In this case, since the reflected wave will have an amplitude 
almost equal to the incident wave, the transmitted wave could have an amplitude almost twice 
as big as that of the incident wave. This explains why the pressure can increase or decrease when
crossing an interface. As a consequence of the different behaviors in the interface, even for 
the one-dimensional case, a complex behavior is observed. These numerical simulations provide 
accurate insight in situations where simple intuition might be insufficient.  

In Figure \ref{fig:dataTBI}, we show from experimental data the shock wave profile in the air before hitting the transwell; however, we are interested in the shape and amplitude of the shock wave in the water just before hitting the endothelial cells sample. In order to do so, we first need to know how relevant the plastic interface is in our model, or if it is relevant at all. Computationally, the plastic interface is hard to model because the width of the plastic is very small ($~mm$) in comparison 
to the shock wave wavelength (meters). The following experiment explores how the width of the plastic interface affects the shock wave profile. Additionally, we show an accurate model can be obtained even when completely ignoring the plastic interface.

\begin{figure*}[h]
  \centering
 \includegraphics[width=0.9\textwidth]{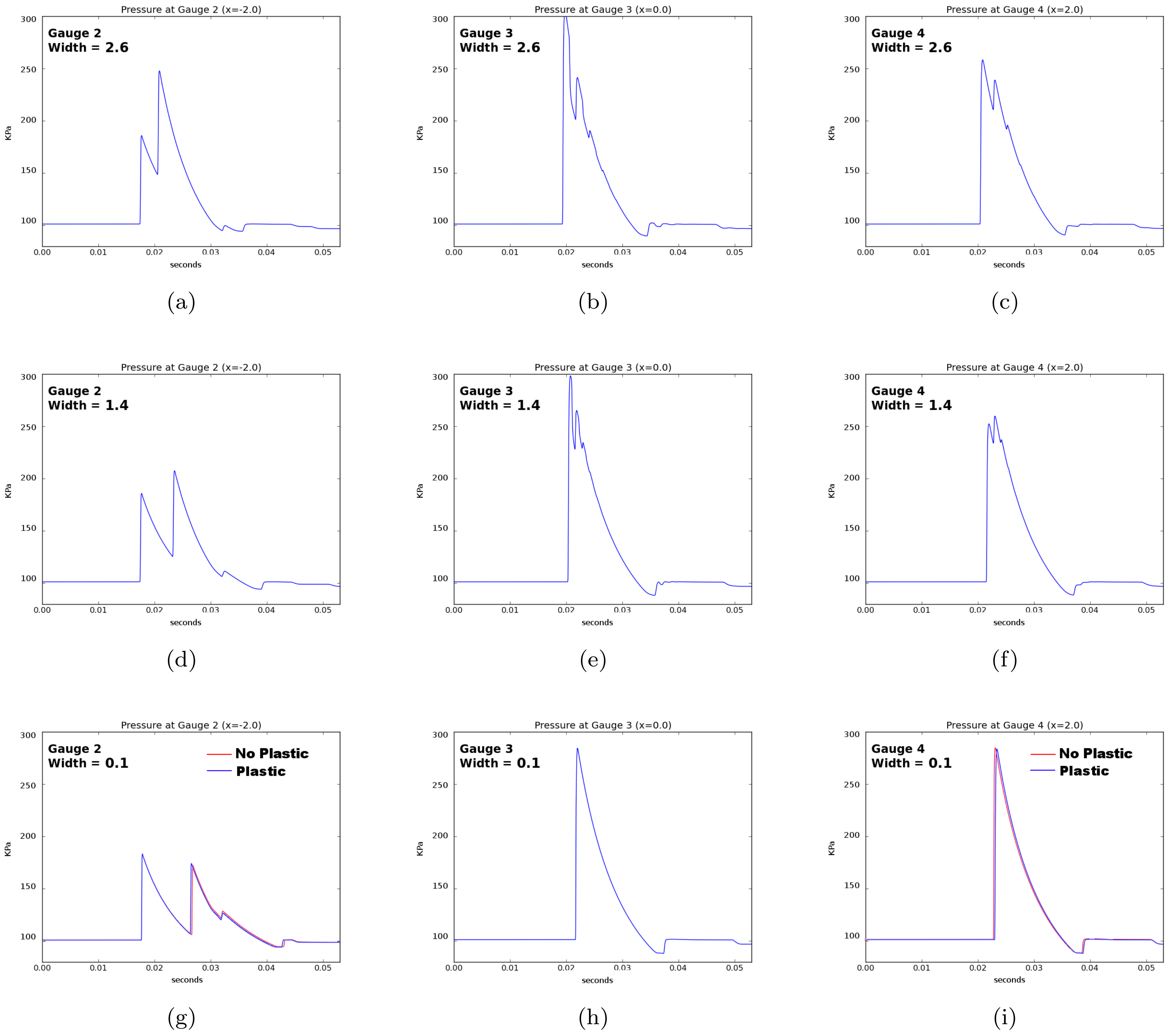} 
  \caption{Pressure (KPa) gauge plots as a function of time (seconds). Each row of figures shows the 
  three gauge plots for three different widths ($2.6$m, $1.4$m and $0.1$m) of the plastic interface, 
  as shown in Table \ref{tab:widthamp}. The plots (g) and (i) for gauge 2 and 4 also show the pressure gauge 
  plots when there is no plastic interface at all; the difference is almost
unnoticeable. }
  \label{fig:TBIgauges}
\end{figure*}

In Figure \ref{fig:TBIanim}, a set of four gauges is shown that allow to measure the 
pressure at certain points as a function of time. The maximum amplitude of the pressure profile was 
measured at gauges $2,3$ and $4$ for different widths of the plastic
interface. The plastic is always assumed to be centered at $x=0$. The results are presented 
in Table \ref{tab:widthamp}. In Figure \ref{fig:TBIgauges}, the full pressure profiles as a 
function of time are shown at the three gauges for three of the plastic widths shown in Table \ref{tab:widthamp}. 

% width = 2*pwidth
\begin{table*}
  \centering
  \begin{tabular}{| l || c | c | c | c |}
    \hline
     Width (m) & Initial(KPa)  & Gauge 2 (KPa) & Gauge 3 (KPa) & Gauge 4 (KPa) \\ \hline
    2.6        & 184.06    & 247.76    & 305.88    & 258.24    \\ \hline
    1.4        & 184.06    & 207.53    & 298.19    & 259.71    \\ \hline
    0.6        & 184.06    & 187.72    & 283.72    & 274.90    \\ \hline
    0.2        & 184.06    & 183.31    & 282.34    & 280.18    \\ \hline  
    0.1        & 184.06    & 183.31    & 284.29    & 283.55    \\
%     2.6        & 184062.09    & 247763.89    & 305880.67    & 258237.71    \\ \hline
%     1.4        & 184062.09    & 207530.14    & 298185.04    & 259705.10    \\ \hline
%     0.6        & 184062.09    & 187724.66    & 283715.90    & 274902.46    \\ \hline
%     0.2        & 184062.09    & 183313.72    & 282343.53    & 280182.40    \\ \hline  
%     0.1        & 184062.09    & 183313.69    & 284287.33    & 283548.49    \\
    \hline
  \end{tabular}
  \caption{The maximum amplitude measured at three pressure gauges for different 
  widths of the plastic interface. The initial shock wave is the same for all cases, and the 
  gauge plots are placed before, inside and after the plastic interface as shown in Figure \ref{fig:TBIanim}.}
  \label{tab:widthamp}
\end{table*}

The results in Table \ref{tab:widthamp} and Figure \ref{fig:TBIgauges} show the maximum amplitude at gauge 2 is reduced 
as the plastic width is decreased. Not surprisingly, this is a consequence of having less interference with the reflected shock wave, since the gauge is farther away from the interface as the plastic width is reduced. This effect is clearly shown on Figures \ref{fig:TBIgauges}a, \ref{fig:TBIgauges}d, \ref{fig:TBIgauges}g. The maximum amplitude at gauge 3 is somewhat diminished at first; however, it seems to be reaching a plateau around $280.0 \text{KPa}$. The behavior at gauge 3 is not trivial; the shock wave bounces back and forth several times, interfering with itself constantly. In Figures \ref{fig:TBIgauges}b, \ref{fig:TBIgauges}e, \ref{fig:TBIgauges}h, we can see the interference becomes so fast that the pressure profile in the plastic seems to converge to a shock wave shape as the plastic width is reduced. At gauge 4, we can observe the interference 
between the set of transmitted shock waves generated by the back and forth reflections within the plastic interface. As the plastic width is reduced, the time elapsed between the transmitted shock waves is reduced and the interference increased. Nonetheless, when the plastic width is very small, the interference becomes so fast that the pressure profile seems to converge again to a shock wave shape, as shown in Figures \ref{fig:TBIgauges}c, \ref{fig:TBIgauges}f, \ref{fig:TBIgauges}i. Furthermore, note the difference in the shock wave shape in Figures \ref{fig:TBIgauges}h, \ref{fig:TBIgauges}i is almost unnoticeable. 
It almost seems like the shock wave is only crossing one interface instead of two. This motivates the next experiment.

\subsection{Air-water interface}
\begin{figure}[H]
  \centering
  \includegraphics[width=0.7\textwidth]{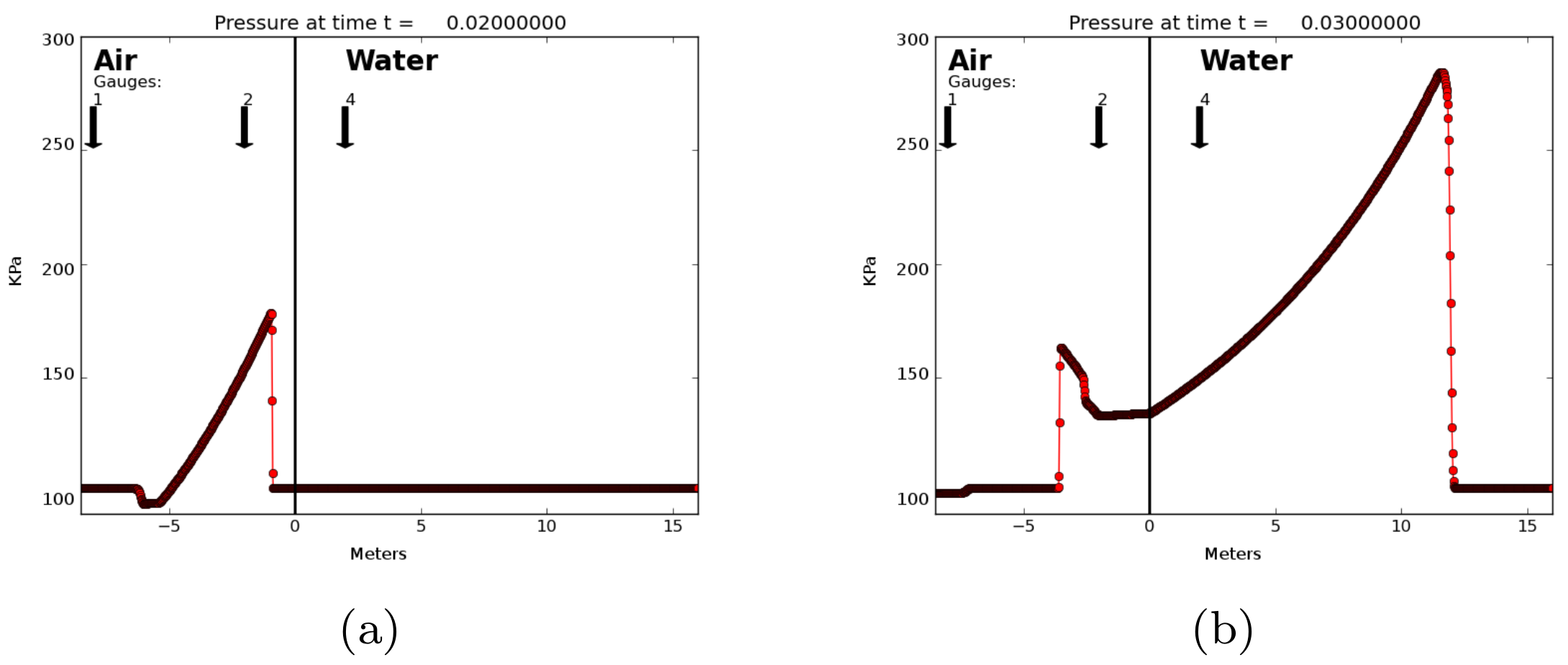} 
  \caption{Shockwave before and after crossing the air-water interface. The arrows indicate the position
of the gauges that measure the pressure as a function of time. The gauges are the same as in Figure \ref{fig:TBIanim}; however, gauge 3 is not shown for plotting aesthetics.}
  \label{fig:TBIanimB}
\end{figure}

In reality, the plastic is so thin that is really unnoticeable on larger scales. 
Furthermore, as the plastic is almost an incompressible medium, one should expect 
it would transfer the shock wave infinitely fast without energy loss. Therefore, 
instead of the triple material interface, let's now consider only an air-water interface. 
The results of this simulation are shown in Figure \ref{fig:TBIanimB}. The gauge plots for 
gauge 2 and 4 are shown in Figures \ref{fig:TBIgauges}g and \ref{fig:TBIgauges}i, along with 
the thin plastic results. The maximum amplitude in each of these gauges is presented in Table \ref{tab:amp}.

\begin{table}[H]
  \centering
  \begin{tabular}{| c | c | c | c |}
    \hline
     Initial (KPa)      & Gauge 2 (KPa) & Gauge 4 (KPa) \\ \hline
     184.06         & 184.40    & 284.26    \\ \hline
%      184062.09         & 183.398.80    & 284258.47    \\ \hline
  \end{tabular}
  \caption{The maximum amplitude measured by two pressure gauges before and after the air-water interface. The initial shock wave
is the same as in the air-plastic-water case.}
  \label{tab:amp}
\end{table}

Comparing the air-water interface results against the ones for the 
smallest plastic width in the air-plastic-water interface case, 
we can observe the percentage error in the maximum pressure amplitude of 
gauge 4 is of $0.38\%$. This is also obvious from the thin plastic and no plastic comparison
in Figures (\ref{fig:TBIgauges}g, \ref{fig:TBIgauges}i). 
An analytic result of this effect for linear acoustics is shown in Section \ref{sec:exactac}. 
These results will allow us to simplify higher dimensional air-plastic-water interface 
problem to a simpler air-water interface. However,
it should be noted that these results for one dimensional interfaces might not be 
as accurate when employed in higher dimensions.

\section{Discussion}
\label{sec:TBIdisc}
Using computational experiments, we showed that there is not a significant difference 
in the transmitted shock wave between the air-plastic-water and the air-water interface. 
Therefore, thin plastic interfaces can be neglected in future computations. 
We also observed an amplification and elongation of the shock wave when crossing the interface. 
This effect is due to the different sound speed and compressibility of the materials. 
The amplitude of the initial pressure wave in air increased by $54 \%$ when measured 
in the water. Since the shock intensity is much higher in water, where the endothelial 
cells sample is, the amplification effect is highly relevant to study injury mechanisms. 
Finally, the amplification effect also occurs when passing from air to plastic or any 
solid material, like bone. This means that a pressure shock wave perceived as weak 
outside in air might be much more intense when traveling through the skull and the 
brain. Future work will involve an extension to three dimensions and implementation of 
two-phase models like those of Pelanti \& Saurel \cite{pelanti2014mixture,saurel2009simple}.
% \cite{Ohl2006}

\bibliographystyle{siam}
\bibliography{extrarefs.bib,library.bib,biorefs.bib}

\end{document}